\documentclass[a4paper,12pt]{article}
\usepackage[english]{babel}
\usepackage[T2A]{fontenc}
\usepackage[cp1251]{inputenc}
\usepackage{amsthm}
\usepackage[tbtags]{amsmath}
\usepackage{amsfonts,amssymb}
\begin{document}

\begin{center}
\textbf{Ekaterina Kompantseva,\footnote{Moscow State Pedagogical University, e-mail: e.kompantseva@gmail.com} Askar Tuganbaev, \footnote{National Research University MPEI; Lomonosov Moscow State University, e-mail: tuganbaev@gmail.com}}

\textbf{Relationships between Almost Completely Decomposable Abelian Groups and Their Multiplication Groups} 
\end{center}

\textbf{Abstract.} For an Abelian group $G$, any homomorphism $\mu\colon G\otimes G\rightarrow G$ is called a \textsf{multiplication} on $G$. The set $\text{Mult}\,G$ of all multiplications on an Abelian group $G$ is an Abelian group with respect to addition. An Abelian group $G$ with multiplication, defined on it, is called a \textsf{ring on the group} 
$G$. Let $\mathcal{A}_0$ be the class of Abelian block-rigid almost completely decomposable groups of ring type with cyclic regulator quotient. 
In the paper, we study relationships between the above groups and their multiplication groups. 
It is proved that groups from $\mathcal{A}_0$ are definable by their multiplication groups. For a rigid group $G\in\mathcal{A}_0$, the isomorphism problem is solved: we describe multiplications from $\text{Mult}\,G$ that define isomorphic rings on $G$. 
We describe Abelian groups that are realized as the multiplication group of some group in $\mathcal{A}_0$. We also describe groups in $\mathcal{A}_0$ that are isomorphic to their multiplication groups.


\textbf{MSC2020 datebase:} 20K30, 20K99, 16B99 

\section{Introduction}

A \textsf{multiplication} on an Abelian group $G$ is a homomorphism $\mu\colon G\otimes G\rightarrow G$. On an Abelian group $G$, the set of all multiplications is an Abelian group with respect to addition; it is denoted by $\text{Mult}\,G$. An Abelian group $G$ with multiplication defined on it is called a \textsf{ring on} $G$. 
The problem of studying the relationships between the structure of an Abelian group and the properties of ring structures on it is very multifaceted and has a long history in algebra; see \cite{BeaP61}, \cite{Fei97}, \cite{Fei00}, \cite{FinR16}, \cite{Gar74}, \cite{Jac82}, \cite{Kom10}, \cite{Kom14}, \cite{Kom15}. 

In this paper, we consider only additively written Abelian groups and the word <<group>> everywhere means <<an Abelian group>>.

The paper is devoted to the study of interrelations between almost completely decomposable Abelian groups and their multiplication groups.

A torsion-free group $G$ of finite rank is called an \textsf{almost completely decomposable} group (\textsf{$ACD$-group}) if $G$ contains a completely decomposable subgroup of finite index. These groups have been studied extensively during the last 50 years (e.g., see \cite{Bla04}, \cite{Bla08}, \cite{BlaM94}, \cite{Bur84}, \cite{DugO93}, \cite{Kom09}, \cite{Lad74}, \cite{Mad00} and other papers). There is a lot of information about these groups, but questions remain.
The book \cite{Mad00} reflects the development of the theory of $ACD$-groups at the time of writing the book.

Any $ACD$-group $G$ contains a special uniquely defined completely decomposable subgroup $\text{Reg}\,G$ of finite index which is a fully invariant subgroup in $G$; it is called the \textsf{regulator} of the group $G$. The regulator of an $ACD$-group can be defined as the intersection of all its completely decomposable subgroups of least index \cite{Bur84}.
The quotient group $G/\text{Reg}\,G$ is called the \textsf{regulator quotient} of the group $G$, and the index of the subgroup $\text{Reg}\,G$ in $G$ is called the \textsf{regulator index}, it is denoted by $n(G)$. 
$ACD$-groups with cyclic regulator quotient are often called \textsf{$CRQ$-groups}.

Let $G$ be an almost completely decomposable group. Then the group $\text{Reg}\,G$ is uniquely, up to isomorphism, representable as a direct sum of torsion-free groups of rank 1 \cite[Proposition 86.1]{Fuc73}. For every type $\tau$, we denote by $\text{Reg}_{\tau}\,G$ the sum of direct summands of rank 1 and of the type $\tau$ in the given direct decomposition of the group $\text{Reg}\,G$. 
The type set
$$
T(G)=T(\text{Reg}\,G)=\{\tau\,|\,\text{Reg}_{\tau}\,G\ne 0\}
$$
is called the \textsf{set of critical types} of the groups $G$ and $\text{Reg}\,G$. 
If $T(G)$ consists of pairwise non-comparable types, then the groups $G$ and $\text{Reg}\,G$ are called \textsf{block-rigid} groups. Moreover, if for any $\tau\in T(G)$, the group $\text{Reg}_{\tau}G$ is of rank $1$, then $G$ and $\text{Reg}\,G$ are called \textsf{rigid} groups. If all types in $T(G)$ are idempotent types then $G$ is called a group \textsf{of ring type}.

In \cite{KomT23}, for a group $G\in\mathcal{A}_0$, we described the group $\text{Mult}\,G$; we proved that it also belongs to the class $\mathcal{A}_0$. In this work, we study interrelations between groups from $\mathcal{A}_0$ and their multiplication groups.
Section 2 is devoted to solving the isomorphism problem for  rings on the rigid group in $\mathcal{A}_0$. For a rigid group $G\in \mathcal{A}_0$, multiplications that define isomorphic rings on $G$ are described (Theorem 2.4). In Section 3, it is proved that groups in the class $\mathcal{A}_0$ are definable by their multiplication groups. 
We say that groups from some class $\mathcal{K}$ \textsf{are definable} by their multiplication groups if for any groups $G,H\in\mathcal{K}$, an isomorphism $\text{Mult}\,G\cong \text{Mult}\,H$ exists if and only if $G\cong H$. We note that in the general, including contemporary, trends in the development of algebra, a significant place is occupied by results that concern the definability of algebraic structures by related structures (see \cite{BlaM23}). In particular, one of the first results in this direction for torsion-free groups was obtained in \cite{BlaIS01}, where it was proved that rigid groups in $\mathcal{A}_0$ are determinable by their endomorphism rings up to near isomorphism, which is a certain weakening of the isomorphism. 
In \cite{BlaIS01}, it is also proved that rigid groups in $\mathcal{A}_0$ are not definable (up to isomorphism) by their endomorphism rings, in general. In the present paper, we show that multiplication groups determine groups from $\mathcal{A}_0$ up to isomorphism  (Theorem 3.5). In Section 4, we describe groups that are realized as the group $\text{Mult}\,G$ for some group $G\in \mathcal{A}_0$ (Theorem 4.2). In particular, it is shown that any rigid group in $\mathcal{A}_0$ can be realized as the multiplication group of some group in $\mathcal{A}_0$ (Corollary 4.3). In addition, we describe groups $G\in \mathcal{A}_0$ which are isomorphic to their multiplication group (Theorem 4.7).

A multiplication $\mu\colon G\otimes G\rightarrow G$ is often denoted by symbols $\times$, $*$, and so on; i.e., $\mu(g_1\otimes g_2)=g_1\times g_2$ for all $g_1,g_2\in G$. Multiplication $\times$ on a group $G$ defines a ring on this group, this ring is denoted by $(G,\times)$.
Let $G$ be a group and $g\in G$. The characteristic and the order of the element $g$ are denoted by $\chi(g)$ and $o(g)$, respectively.
The rank of the group $G$ and the divisible hull of $G$ are denoted by $r(G)$ and $\tilde{G}$, respectively. If $S\subseteq G$, then $\langle S\rangle$ is the subgroup of the group $G$ generated by the set $S$. 
The elements of the direct product $\prod_{i\in I}G_i$ of groups are written as $(g_i)_{i\in I}$, $g_i\in G_i$. If $I_1\subseteq I$, then for simplicity the subgroup $\{(g_i)_{i\in I}\in \prod_{i\in I}G_i\,|\,g_i=0$ for all $i\notin I_1\}$ of the group $\prod_{i\in I}G_i$ is identified with the group $\prod_{i\in I_1}G_i$, and its elements are written as $(g_i)_{i\in I_1}$. Even if the set $I$ is finite, we formally distinguish the groups $\prod_{i\in I}G_i$ and $\oplus_{i\in I}G_i$; in the last case, the elements of the direct sum are denoted by $\sum_{i\in I}g_i$, $g_i\in G_i$.
As usual, $\mathbb{N}$, $\mathbb{P}$ are the sets of all positive integers and all prime numbers, respectively, 
$\mathbb{Z}$ is the group (the ring) of integers, $\mathbb{Q}$ is the group (the field) of rational numbers. If $R$ is a ring with identity element, then $Re$ is the cyclic module over $R$ generated by the element $e$. If $S$ is a finite subset in $\mathbb{Z}$, then $\text{gcd}(S)$ is the greatest common divisor of all numbers from $S$, $\text{lcm}(S)$ is the least common multiple of numbers from $S$.
If $P_1\subseteq \mathbb{P}$, then a \textsf{$P_1$-number} is a non-zero integer whose every prime divisor (if it exists) is contained in $P_1$,
a \textsf{$P_1$-fraction} is a rational number that can be represented as a fraction whose numerator and denominator are $P_1$-numbers.
It follow from the definition that $1$ is a $P_1$-number for every $P_1\subseteq \mathbb{P}$.
For any type $\tau$, we denote by
$$
P_{\infty}(\tau)=\{p\in \mathbb{P}\,|\, \tau(p)=\infty\}, \quad 
P_0(\tau)=\mathbb{P}\setminus P_{\infty}(\tau).
$$
For all definitions and notation, unless otherwise stated, we refer to the books \cite{Fuc73}, \cite{Fuc15} and \cite{KMT03}.

\section{Isomorphism Proplem}\label{section2}
 
This section is devoted to the isomorphism problem, which consists in describing multiplications on the group $G\in \mathcal{A}_0$ that define isomorphic rings on $G$.

Let $G$ be a reduced block-rigid $CRQ$-group of ring type with regulator $A$ and regulator quotient $G/A=\langle d+A\rangle$, i.e. $G=\langle d,A\rangle$. Let $n=n(G)$ be the regulator index of the group $G$ and let $T(G)$ be the set of critical types.

Denote $\text{Reg}_{\tau}\,G=A_{\tau}$, then the group $A$ can be represented in the form $A=\oplus_{\tau\in T(G)}A_{\tau}$. According to \cite[Proposition 2.4.11]{Mad00}, such decomposition of a completely decomposable group is unique if and only if $A$ is a block-rigid group. For divisible hulls $\tilde G$, $\tilde A$ and $\tilde A_{\tau}$ of the groups $G$, $A$ and $A_{\tau}$, respectively, we have equalities 
$$
\tilde G=\tilde A=\oplus_{\tau\in T(G)}\tilde A_{\tau}.
$$
For $\tau\in T(G)$, we denote by $\pi_{\tau}$ the projection of the group $\tilde G$ onto $\tilde A_{\tau}$. 

For a group $G\in \mathcal{A}_0$, positive integers $m_{\tau}=m_{\tau}(G)$ ($\tau\in T(G)$) are defined as follows, \cite{DugO93}. 
Let $d_{\tau}=\pi_{\tau}(d)\in \tilde{A}_{\tau}$ and let $m_{\tau}=o(d_{\tau}+A)$ be the order of the element $d_{\tau}+A$ in the torsion group $\tilde{A}/A$. We note that the numbers $m_{\tau}$ do not depend on the choice of the element $d$ and
$n(G)=o(d+A)=\text{lcm}\{m_{\tau}\,|\,\tau\in T(G)\}$.
In addition, according to \cite{DugO93}, the set $\{m_{\tau}\,|\,\tau\in T(G)\}$ is a system of near isomorphism invariants of the group $G$.

In \cite[Theorem 3.5]{BlaM94}, it is proved that for any group $G\in \mathcal{A}_0$, there exists a direct decomposition 
$$
G=G_1\oplus C,\eqno (2.1)
$$

\noindent where $C$ is a block-rigid completely decomposable group and $G_1$ is a rigid $CRQ$-group which satisfies the following conditions:
$$
\tau\in T(G_1) \text{ if and only if } m_{\tau}(G)> 1,\eqno (2.1')
$$
$$
m_{\tau}(G_1)=m_{\tau}(G) \text{ for all } \tau\in T(G_1).\eqno (2.1'')
$$

Decomposition $(2.1)$, which satisfies conditions $(2.1')$ and $(2.1'')$, is called a \textsf{main decomposition} of the group $G$. The group $G_1$ in a main decomposition of the group $G$ does not contain completely decomposable direct summands; such groups are called \textsf{clipped}. We note that a main decomposition of a $CRQ$-group is not uniquely defined \cite{BlaM94}. 
We assume that the main decomposition of the group $G$ is fixed everywhere.

We set $T_0(G)=\{\tau\in T(G)\,|\, m_{\tau}>1\}$. Then $T_0(G)$ is the set of critical types of the clipped direct summand in any main decomposition of the group $G$.
We set $D=\{d\in G_1\,|\,G/A=\langle d+A\rangle\}$. We denote by $B$ the regulator of the group $G_1$. Then $T(G_1)=T(B)=T_0(G)$ and $\tilde{G}_1=\tilde{B}$. Let $d\in D$; then there exists a system $E_0=\{e_0^{(\tau)}\in B_{\tau}\,|\, \tau\in T_0(G)\}$ such that
$$
B=\oplus_{\tau\in T_0(G)}R_{\tau}e_0^{(\tau)}, \eqno (2.2)
$$
where $R_{\tau}$ is a unitary subring of the field $\mathbb{Q}$, the type of the additive group of $R_{\tau}$ is equal to $\tau$, and an element $d$ in the group $\tilde{B}$ can be represented in the form
$$
d=\sum_{\tau\in T_0(G)}\dfrac{s_{\tau}}{m_{\tau}}e_0^{(\tau)}, \eqno (2.3)
$$
where integers $m_{\tau}$ and $s_{\tau}$ satisfy the following conditions: 
$$
\text{gcd}(s_{\tau},m_{\tau})=1 \text{ for all } \tau\in T_0(G).
\eqno (2.3')
$$
$$
s_{\tau} \text{ and } m_{\tau} \text{ are }P _0(\tau)\text{-numbers for every } \tau\in T_0(G).\eqno(2.3'')
$$

\noindent A system $E_0=\{e_0^{(\tau)}\in B_{\tau} \,|\, \tau\in T_0(G)\}$, which satisfies conditions $(2.2)$ and $(2.3)$, is called an \textsf{$rc$-basis} of the group $G$ defined by the element $d$. We note that the pair $(d,E_0)$ uniquely defines numbers $s_{\tau}$ ($\tau\in T_0(G)$).
For the block-rigid $CRQ$-group $G$, the equality $(2.3)$ is called the \textsf{standard representation} of $G$ related to the pair $(d,E_0)$. 

For $\tau\in T(G)$, we denote by $k_{\tau}=r(A_{\tau})$; we set 
$$
I_{\tau}=\begin{cases}
\{0,1,\ldots,k_{\tau}-1\} \text{ for } \tau\in T_0(G)\\
\{1,2,\ldots,k_{\tau}\} \text{ for } \tau\notin T_0(G).
\end{cases}
$$
Then there exists a system $E=\{e_i^{(\tau)}\in A_{\tau}\,|\,i\in I_{\tau}, \;\tau\in T(G)\}$ such that 
$$
A=\oplus_{\tau\in T(G)}\oplus_{i\in I_{\tau}}R_{\tau}e_i^{(\tau)},
$$
where $R_{\tau}$ is a unitary subring of the field $\mathbb{Q}$, the type of the additive group of $R_{\tau}$ is equal to $\tau$. This  system $E$ is called an \textsf{$r$-basis} of the group $G$ if its subsystem $E_0=\{e_0^{(\tau)}\,|\,\tau\in T_0(G)\}$ is an $rc$-basis of the group $G$.

For an $rc$-basis $E_0$, we denote
$$
D(E_0)=\{d\in D\,|\, rc\text{-basis } E_0 \text{ is determined by the element } d\}. 
$$ 
Then $D(E_0)\ne \varnothing$ and $D(E_0)$ may contain more than one element (see \cite[Remark 2.1]{KomT23}).

To describe the group $\text{Mult}\,G$, we define the following groups. 
For any $\tau\in T(G)$, we set\\
$M_{\tau}^{(0)}=M_{k_{\tau}}(A_{\tau})$ is the additive group of square matrices of order $k_{\tau}$ with elements in $A_{\tau}$, we also set 
$$
M_{\tau}^{(1)}=\begin{bmatrix}
m_{\tau}A_{\tau}&m_{\tau}A_{\tau}&\ldots&m_{\tau}A_{\tau}\\
m_{\tau}A_{\tau}&A_{\tau}&\ldots&A_{\tau}\\
\ldots&\ldots&\ldots&\ldots\\
m_{\tau}A_{\tau}&A_{\tau}&\ldots&A_{\tau}
\end{bmatrix}\subseteq M_{\tau}^{(0)},
$$
where the notation $[\ldots]$ means the set of matrices of a certain form, $m_{\tau}=m_{\tau}(G)$,
$$
M_{\tau}^{(2)}=\begin{bmatrix}
m_{\tau}^2A_{\tau}&m_{\tau}A_{\tau}&\ldots&m_{\tau}A_{\tau}\\
m_{\tau}A_{\tau}&A_{\tau}&\ldots&A_{\tau}\\
\ldots&\ldots&\ldots&\ldots\\
m_{\tau}A_{\tau}&A_{\tau}&\ldots&A_{\tau}
\end{bmatrix}
\subseteq M_{\tau}^{(1)}.
$$
We set 
$$
M^{(0)}=\prod_{\tau\in T(G)}M_{\tau}^{(0)},\;
M^{(1)}=\prod_{\tau\in T(G)}M_{\tau}^{(1)},\;
M^{(2)}=\prod_{\tau\in T(G)}M_{\tau}^{(2)}.
$$
Then $M^{(2)}\subseteq M^{(1)}\subseteq M^{(0)}$.

For the standard representation $(2.3)$ of $G$ related to the pair $(d,E_0)$ and for every $\tau\in T(G)$, we consider the matrix
$$
X^{(\tau)}=X^{(\tau)}(d,E_0)=\begin{pmatrix}
m_{\tau}s_{\tau}^{-1}e_0^{(\tau)}&0&\ldots&0\\
0&0&\ldots&0\\
\ldots&\ldots&\ldots&\ldots\\
0&0&\ldots&0
\end{pmatrix}\in M_{\tau}^{(1)},\;\text{if } \tau\in T_0(G),
$$
where $s_{\tau}^{-1}$ is an integer that is inverse to $s_{\tau}$ modulo $m_{\tau}$,
$$
X^{(\tau)}=\begin{pmatrix}
0&0&\ldots&0\\
0&0&\ldots&0\\
\ldots&\ldots&\ldots&\ldots\\
0&0&\ldots&0
\end{pmatrix}\in M_{\tau}^{(1)},\;\text{if } \tau\notin T_0(G).
$$
We set 
$$
X=X(d,E_0)=\left(X^{(\tau)}\right)_{\tau\in T(G)}=\left(X^{(\tau)}\right)_{\tau\in T_0(G)}\in M^{(1)},
$$
$$
M(d,E_0)=\langle X,M^{(2)}\rangle\subseteq M^{(1)}.
$$
We note that the integral solutions of the congruence $s_{\tau}x\equiv 1\,(\text{mod}\,m_{\tau})$ form a class of residues modulo $m_{\tau}$. Therefore, the set $M(d,E_0)$ does not depend on the choice of numbers $s_{\tau}^{-1}$ when defining $X$.
Also note that if $\tau\in T(C)\setminus T_0(G)$, then $m_{\tau}=1$ by $(2.1')$. Therefore, we have
$$
M_{\tau}^{(2)}=M_{\tau}^{(1)}=M_{\tau}^{(0)}
$$
in this case.

\textbf{Remark 2.1.} In \cite{KomT23}, it is proved that the group $M(d,E_0)$ does not depend on the choice of the element $d\in D(E_0)$. Therefore, we denote $M(d,E_0)=M_G(E_0)$. In \cite[Theorem 2.8]{KomT23}, it is proved that if $G\in \mathcal{A}_0$ and $E_0$ is an $rc$-basis of the group $G$, then $\text{Mult}\,G\cong M_G(E_0)$.~$\square$

We want to solve the isomorphism problem for rigid groups in the class $\mathcal{A}_0$. Now it can be formulated as follows: 
Describe the elements of $M_G(E_0)$ that define isomorphic rings on $G$.

In this section, later $G$ is a rigid group in $\mathcal{A}_0$, $T=T(G),$ $m_{\tau}=m_{\tau}(G),$ $T_0=T_0(G)$. 
Then an $r$-basis of the group $G$ can be written in the form $E=\{e_{\tau}\,|\,\tau\in T\}$. The regulator of $G$ can be written in the form $A=\text{Reg}\,G=\oplus_{\tau\in T}R_{\tau}e_{\tau}$, where $R_{\tau}$ is a unitary subring of the field $\mathbb{Q}$, the type of the additive group of $R_{\tau}$ is equal to $\tau$. In addition, the system $E_0=\{e_{\tau}\,|\,\tau\in T_0\}$ is an $rc$-basis of the group $G$. Let a standard representation of the group $G$ be of the form 
$$
d=\sum_{\tau\in T_0}\dfrac{s_{\tau}}{m_{\tau}}e_{\tau}.\eqno (2.4)
$$
According to Remark 2.1, the group $M_G(E_0)\cong \text{Mult}\,G$ is of the form
$$
M_G(E_0)=\langle X,M^{(2)}\rangle\subseteq \prod_{\tau\in T}m_{\tau}R_{\tau}e_{\tau},\eqno (2.5)
$$ 
where 
$$
X=(m_{\tau}s_{\tau}^{-1}e_{\tau})_{\tau\in T_0}\in\prod_{\tau\in T}m_{\tau}R_{\tau}e_{\tau},\; M^{(2)}=\prod_{\tau\in T}m_{\tau}^2R_{\tau}e_{\tau},
$$
$s_{\tau}^{-1}$ is an integer that is inverse to $s_{\tau}$ modulo $m_{\tau}$; we recall that $m_{\tau}=1$ for $\tau\in T\setminus T_0$.
In addition, the isomorphism $\text{Mult}\,G\cong M_G(E_0)$ takes each multiplication $\times\in\text{Mult}\,G$ to an element $U_{\times}=(u_{\tau}e_{\tau})_{\tau\in T}\in M_G(E_0)$ with $e_{\tau}\times e_{\tau}=u_{\tau}e_{\tau}$ for $\tau\in T$.

We consider the ring $R_G=\prod_{\tau\in T}R_{\tau}$. Then $R_G$ is a unital ring. An element $c=(c_{\tau})_{\tau\in T}\in R_G$ is called \textsf{$P_{\infty}$-element} if $c_{\tau}$ is a $P_{\infty}(\tau)$-fraction for every $\tau\in T$. We denote by $1_{\tau}$ the identity element of the ring $R_{\tau}$, $1_{T_0}=(1_{\tau})_{\tau\in T_0}\in R_G$. In $R_G$, we consider the subring 
$$
K_G=\mathbb{Z}1_{T_0}+\prod_{\tau\in T}m_{\tau}R_{\tau}.
$$
Then each of the groups $\prod_{\tau\in T}\mathbb{Q}e_{\tau}$ and $\oplus_{\tau\in T}\mathbb{Q}e_{\tau}$ is a module over $R_G$ (and, consequently, over $K_G$) if for any 
$$
r=(r_{\tau})_{\tau\in T}\in R_G,\; a=(ae_{\tau})_{\tau\in T}\in\prod_{\tau\in T}\mathbb{Q}_{\tau}e_{\tau},\; b=\sum_{\tau\in T}b_{\tau}e_{\tau}\in \oplus_{\tau\in T}\mathbb{Q}e_{\tau},
$$
we set $ra=(r_{\tau}a_{\tau}e_{\tau})_{\tau\in T}$ and $rb=\sum_{\tau\in T}r_{\tau}b_{\tau}e_{\tau}$.

For an arbitrary unital ring $R$, we denote by $R^*$ its multiplicative group.

\textbf{Lemma 2.2.} Let $G$ be a rigid group from $\mathcal{A}_0$, $T(G)=T$, and let $E=\{e_{\tau}\,|\,\tau\in T\}$ be a $r$-basis of the group $G$. Then the following assertions are true.

\textbf{1.} $G$ is a submodule of the module $\tilde{G}=\oplus_{\tau\in T}\mathbb{Q}e_{\tau}$ over $K_G$.

\textbf{2.} $M_G(E_0)$ is a submodule of the module $\prod_{\tau\in T}\mathbb{Q}e_{\tau}$ over $K_G$.

\textbf{3.} The multiplicative group $R_G^*$ coincides with the set of all $P_{\infty}$-elements of the ring $R_G$.

\textbf{4.} For the multiplicative group $K_G^*$ of the ring $K_G$, we have $K_G^*=K_G\cap R_G^*$.

\textbf{Proof.} \textbf{1.} Let the group $G$ have the standard representation $(2.4)$. It is easy to see that $R_GA\subseteq A\subseteq G$. Let $c=\gamma 1_{T_0}+(m_{\tau}y_{\tau})_{\tau\in T}\in K_G$, where $\gamma\in\mathbb{Z}$ and $y_{\tau}\in R_{\tau}$. Then 
$$
cd=\gamma d+\sum_{\tau\in T_0}\left(m_{\tau}y_{\tau}\dfrac{s_{\tau}}{m_{\tau}}\right)e_{\tau} =\gamma d+ \sum_{\tau\in T_0}y_{\tau}s_{\tau}e_{\tau}\in\langle d,A\rangle=G.
$$
Consequently, $G$ is a submodule of the module $\oplus_{\tau\in T}\mathbb{Q}e_{\tau}$ over $K_G$.

\textbf{2.} Let the group $M_G(E_0)$ be represented in the form $(2.5)$. It is easy to see that $R_GM^{(2)}\subseteq M^{(2)}\subseteq M_G(E_0)$. Let
$c=\gamma 1_{T_0}+(m_{\tau}y_{\tau})_{\tau\in T}\in K_G$, where $\gamma\in \mathbb{Z}$, $y_{\tau}\in R_{\tau}$. Then
$$
cX=\gamma X+(m_{\tau}^2y_{\tau}s_{\tau}^{-1}e_{\tau})_{\tau\in T_0}\in\langle X,M^{(2)}\rangle=M_G(E_0).
$$
Consequently, $M_G(E_0)$ is a submodule of the module $\prod_{\tau\in T}\mathbb{Q} e_{\tau}^*$.

\textbf{3.} The assertion follows from the equality $R_G^*=\prod_{\tau\in T}R_{\tau}^*$.

\textbf{4.} The inclusion $K_G^*\subseteq K_G\cap R_G^*$ is obvious. 
Let us prove the reverse inclusion. Let $c\in K_G\cap R_G^*$, 
$c=(c_{\tau})_{\tau\in T}=(\gamma+m_{\tau}y_{\tau})_{\tau\in T_0}+
(c_{\tau})_{\tau\in T\setminus T_0}$, where $\gamma\in\mathbb{Z}$, $y_{\tau}\in R_{\tau}$ for $\tau\in T_0$, $c_{\tau}$ is a $P_{\infty}(\tau)\text{-fraction}$ for $\tau\in T$.
Then for any $\tau\in T$, there exists a $c_{\tau}^{-1}\in R_{\tau}$.

Let $\tau\in T_0$, $y_{\tau}=\dfrac{y_1}{y_2}$, where $y_1\in\mathbb{Z}$, $y_2$ is a $P_{\infty}(\tau)$-number; then $c_{\tau}=\gamma+m_{\tau}y_{\tau}=\dfrac{\gamma y_2+m_{\tau}y_1}{y_2}$. Therefore, $\gamma y_2+m_{\tau}y_1$ is a $P_{\infty}(\tau)$-number; hence, $\text{gcd}(\gamma y_2+m_{\tau}y_1,m_{\tau})=1$. 
Consequently, $\text{gcd}(\gamma,m_{\tau})=1$ for every $\tau\in T$; therefore $\gamma$ is co-prime to $n=n(G)$. Let $\gamma^{-1}$ be an integer that is inverse to $\gamma$ modulo $n$; then $\gamma\gamma^{-1}=1+m_{\tau}v_{\tau}$ for some $v_{\tau}\in\mathbb{Z}$. It is directly verified that $c_{\tau}^{-1}=\gamma^{-1}-m_{\tau}z_{\tau}$, where $z_{\tau}=c_{\tau}^{-1}(\gamma^{-1}y_{\tau}+v_{\tau})\in R_{\tau}$. Consequently, 
$$
c^{-1}=\gamma^{-1}1_{T_0}-(m_{\tau}z_{\tau})_{\tau\in T_0}+(c_{\tau}^{-1})_{\tau\in T\setminus T_0}\in K_G;
$$
therefore, $c\in K_G^*$.~$\square$

\textbf{Lemma 2.3.} Let $G$ be a rigid group in $\mathcal{A}_0$. 

If $f$ is an endomorphism (resp., an automorphism) of the group $G$, then there is an element $c\in K_G$ (resp., $c\in K_G^*$) such that $f(x)=cx$ for all $x\in G$.

If $c\in K_G$ (resp., $c\in K_G^*$), then there is an endomorphism (resp., an automorphism) $f\colon G\to G$ such that $f(x)=cx$ for all $x\in G$.

\textbf{Proof.} Let $\{e_{\tau}\,|\,\tau\in T\}$ be an $r$-basis of the group $G$ and $f$ be an endomorphism (resp., an automorphism) of the group $G$; then for every $\tau\in T$, the element $f(e_{\tau})$ can be represented in the form $f(e_{\tau})=c_{\tau}e_{\tau}$, where $c_{\tau}\in R_{\tau}$, since $A_{\tau}=R_{\tau}e_{\tau}$ is a fully invariant subgroup of the group $G$. We set $c=(c_{\tau})$; then for any $x\in G$ we have $f(x)=cx$. It follow from Theorem 6.2 and Theorem 6.4 in \cite{MadS00} that $c\in K_G$ (resp., $c\in K_G^*$).

Let $c\in K_G$. By Lemma 2.2(1), there is an endomorphism $f\colon G\to G$ such that $f(x)=cx$ for all $x\in G$. 
If $c\in K_G^*$, then an endomorphism $f^{-1}\colon G\to G$ such that $f^{-1}(x)=c^{-1}x $ for all $x\in G$ is inverse to $f$. Therefore, $f$ is an automorphism of the group $G$.~$\square$

\textbf{Theorem 2.4 (The isomorphism theorem).}\\
Let $G$ be a rigid group in $\mathcal{A}_0$ with an $rc$-basis $E_0$. Let elements $U_*,U_{\times}\in M_G(E_0)$ determine the rings $(G,*)$ and $(G,\times)$, respectively. Then $(G,*)\cong (G,\times)$ if and only if there exists an element $c\in K_G^*$ such that $U_*=cU_{\times}$.

\textbf{Proof.} Let $U_*=(u_{\tau}e_{\tau})_{\tau\in T}$, where $u_{\tau}e_{\tau}=e_{\tau}*e_{\tau}$, $u_{\tau}\in R_{\tau}$ $(\tau\in T)$, $U_{\times}=(v_{\tau}e_{\tau})_{\tau\in T}$, where $v_{\tau}e_{\tau}=e_{\tau}\times e_{\tau}$, $v_{\tau}\in R_{\tau}$ $(\tau\in T)$.

If $f\colon (G,*)\to (G,\times)$ is a ring isomorphism, then $f$ is an automorphism of the additive group $G$. It follows from Lemma 2.3 that there exists an element $c\in K_G^*$ such that $f(x)=cx$ for every $x\in G$. Let $\tau\in T$; then 
$$
f(e_{\tau}*e_{\tau})=f(u_{\tau}e_{\tau})=u_{\tau}f(e_{\tau})=c_{\tau}u_{\tau}e_{\tau}.
$$
On the other hand,
$$
f(e_{\tau}*e_{\tau})=f(e_{\tau})\times f(e_{\tau})=c_{\tau}e_{\tau}\times c_{\tau}e_{\tau}=c_{\tau}^2(e_{\tau}\times e_{\tau})=c_{\tau}^2v_{\tau}e_{\tau}.
$$
Consequently, $c_{\tau}u_{\tau}e_{\tau}=c_{\tau}^2v_{\tau}e_{\tau}$. Therefore, $u_{\tau}e_{\tau}=c_{\tau}v_{\tau}e_{\tau}$ for all $\tau\in T$; consequently, $U_*=cU_{\times}$.

Now let $U_*=cU_{\times}$ for some $c=(c_{\tau})_{\tau\in T}\in K_G^*$. Then
$$
u_{\tau}e_{\tau}=c_{\tau}v_{\tau}e_{\tau}\, \text{ for any }\, \tau\in T. \eqno (2.6)
$$
By Lemma 2.3, there is a group automorphism $f\colon G\to G$ such that $f(x)=cx$ for $x\in G$. In addition, $f$ also is an isomorphism of the ring $(G,*)$ onto the ring $(G,\times)$. Indeed, by $(2.6)$ for any $\tau\in T$ we have $f(e_{\tau}*e_{\tau})=$
$$
=f(u_{\tau}e_{\tau})=c_{\tau}u_{\tau}e_{\tau}=c_{\tau}^2v_{\tau}e_{\tau}=c_{\tau}^2(e_{\tau}\times e_{\tau})=
(c_{\tau}e_{\tau})\times (c_{\tau}e_{\tau})= f(e_{\tau})\times f(e_{\tau}).\,\square
$$

\section{Definability of CRQ-Groups by Their Multiplcation Groups}\label{section3}

We present some familiar results about multiplication groups of groups in the class $\mathcal{A}_0$.

\textbf{Theorem 3.1 \cite{KomT23}.} Let $G\in \mathcal{A}_0$.
Then for the group $\text{Mult}\,G$, the following statements hold.

\textbf{1.} The group $\text{Mult}\,G$ is a block-rigid $CRQ$-group of ring type with regulator $\text{Hom}(G\otimes G,A)\cong M^{(2)}$.

\textbf{2.} $T(\text{Mult}\,G)=T(G)$.

\textbf{3.} We have $m_{\tau}(\text{Mult}\,G)=m_{\tau}(G)$ for any $\tau\in T(G)$, and, as a consequence, $T_0(\text{Mult}\,G)=T_0(G)$,~ $n(\text{Mult}\,G)=n(G)$.

\textbf{4.} $r(\text{Reg}_{\tau}(\text{Mult}\,G))=(r(\text{Reg}_{\tau}G))^3$ for any $\tau\in T(G)$.

\textbf{5.} Let $E_0$ be an $rc$-basis of the group $G$. We set
$$
F_0^{(\tau)}=\begin{pmatrix}
m_{\tau}^2e_0^{(\tau)}&0&\ldots&0\\
0&0&\ldots&0\\
\ldots&\ldots&\ldots&\ldots\\
0&0&\ldots&0
\end{pmatrix}\in M_{\tau}^{(2)}.
$$
Then the system $\{F_0^{(\tau)}\,|\,\tau\in T_0(G)\}$ is a one of $rc$-bases of the group $M_G(E_0)\cong \text{Mult}\,G$.

\textbf{6.} If the standard representation of the group $G$ is of the form (2.3), then for every $\tau\in T_0(G)$, there exists a $P_0(\tau)$-number $s_{\tau}^{-1}$ that is inverse to $s_{\tau}$ modulo $m_{\tau}$. 
One of the standard representations of the group $\text{Mult}\,G$ is of the form
$X=\left(\dfrac{s_{\tau}^{-1}}{m_{\tau}}F_0^{(\tau)}\right)_{\tau\in T_0(G)}$.~$\square$

To study relationships between a group in the class $\mathcal{A}_0$ and its multiplication group, we need a generalization of the notion of isomorphism. 
Let $G$ and $H$ be two torsion-free groups of finite rank. According to \cite{Mad00} the groups $G$ and $H$ are called \textsf{near isomorphic} (we designate $G\cong_{nr}H$) if for every prime number $p$, there exists a monomorpism $f_p\colon G\to H$ such that the index of the subgroup $f_p(G)$ in $H$ is finite and it is not divided by $p$.

\textbf{Remark 3.2.}\\
\textbf{1.} According to \cite[Proposition 86.1]{Fuc73}, two completely decomposable groups $G$ and $H$ are isomorphic if and only if $T(G)=T(H)$ and $r(G_{\tau})=r(H_{\tau})$ for all $\tau\in T(G)$.

\textbf{2.} Let $G$ and $H$ be two block-rigid $CRQ$-groups of ring type. According to \cite{BlaM94}, the groups $G$ and $H$ are near isomorphic if and only if their regulators are isomorphic and $m_{\tau}(G)=m_{\tau}(H)$ for all types $\tau\in T(G)=T(H)$.~$\square$

\textbf{Theorem 3.3.} Let $G_1,G_2\in\mathcal{A}_0$. The groups $G_1$ and $G_2$ are near isomorphic if and only if the groups $\text{Mult}\,G_1$ and $\text{Mult}\,G_2$ are near isomorphic.

\textbf{Proof.} 
Let $G_1,G_2\in\mathcal{A}_0$. By Theorem 3.1,
we have each of the following conditions for every $i\in \{1,2\}$:
$$
\begin{cases}
T(\text{Mult}\,G_i)=T(G_i),\\
r[\text{Reg}_{\tau}(\text{Mult}\,G_i)]=(r(\text{Reg}_{\tau}G_i))^3 \text{ for all } \tau\in T(G_i),\\
m_{\tau}(\text{Mult}\,G_i)=m_{\tau}(G_i) \text{ for all } \tau\in T(G_i).
\end{cases} \eqno(3.1)
$$
By Remark 3.2, the groups $G_1$ and $G_2$ are near isomorphic if and only if 
$$
\begin{cases}
T(G_1)=T(G_2),\\
r(\text{Reg}_{\tau}(G_1))=r(\text{Reg}_{\tau}G_2) \text{ for any } \tau\in T(G_i),\\
m_{\tau}(G_1)=m_{\tau}(G_2) \text{ for any } \tau\in T(G_i).
\end{cases} \eqno(3.2)
$$
Similarly, the groups $\text{Mult}\,G_1$ and $\text{Mult}\,G_2$ are near isomorphic if and only if 
$$
\begin{cases}
T(\text{Mult}\,G_1)=T(\text{Mult}\,G_2),\\
r[\text{Reg}_{\tau}(\text{Mult}\,G_1)]=r[(\text{Reg}_{\tau}(\text{Mult}\,G_2)] \text{ for any } \tau\in T(G_i),\\
m_{\tau}(\text{Mult}\,G_1)=m_{\tau}(\text{Mult}\,G_2) \text{ for any } \tau\in T(G_1)
\end{cases} \eqno(3.3)
$$
According to $(3.1)$, conditions $(3.2)$ and $(3.3)$ are equivalent.~$\square$

In order to describe isomorphic groups in the near isomorphism class of groups in the class $\mathcal{A}_0$, we introduce the following notation.
Let $G$ be a group in the class $\mathcal{A}_0$ with regulator $A$, regulator index $n$ and invariants of near isomorphism $m_{\tau}$ ($\tau\in T(G)$). Let
$$
T_0=T_0(G)=\{\tau\in T(G)\,|\,m_{\tau}>1\},
$$
$$
T_1=T_1(G)=\{\tau\in T(G)\,|\,m_{\tau}>1,\,r(A_{\tau})=1\}.
$$
We set 
$$
S=\prod_{\tau\in T_1}\mathbb{Z}_{m_{\tau}}\,-\,\text{the direct product of the rings } \mathbb{Z}_{m_{\tau}}=\mathbb{Z}/m_{\tau}\mathbb{Z},
$$
$$
\overline{1_{\tau}}=1+m_{\tau}\mathbb{Z}\in\mathbb{Z}_{m_{\tau}}.$$
Then $(\overline{1_{\tau}})_{\tau\in T_1}$ is the identity element of the ring $S$. 
We write the elements of the ring $S$ in the form $(x_{\tau}\overline{1_{\tau}})_{\tau\in T_1}$, where $x_{\tau}\in\mathbb{Z}$. 
We consider the multiplicative group $S^*$ of the ring $S$ which is equal to the direct product $\prod_{\tau\in T_1}\mathbb{Z}_{m_{\tau}}^*$ of multiplicative groups $\mathbb{Z}_{m_{\tau}}^*$. Then
$$
S^*=\{(x_{\tau}\overline{1_{\tau}})_{\tau\in T_1}\in S\,|\,
x_{\tau}\in\mathbb{Z},\,\text{gcd}(x_{\tau},m_{\tau})=1 \text{ for }\tau\in T_1\}.
$$
In the group $S^*$, we define two subgroups
$$
\Gamma=\{(\alpha\overline{1_{\tau}})_{\tau\in T_1}(\beta\overline{1_{\tau}})_{\tau\in T_1}^{-1}\,|\,
\alpha,\beta\in\mathbb{Z},\,\text{gcd}(\alpha,n)=\text{gcd}(\beta,n)=1\},
$$
$$
V_{\infty}=\{(x_{\tau}\overline{1_{\tau}})_{\tau\in T_1} (y_{\tau}\overline{1_{\tau}})_{\tau\in T_1}^{-1}\,|\,x_{\tau},y_{\tau} - P_{\infty}(\tau)\text{-numbers} \text{ for }\tau\in T_1\}.
$$
\textbf{Remark 3.4.} Let $G$ and $H$ be two near isomorphic groups from the class $\mathcal{A}_0$. By Remark 3.2, we have
$$
T(G)=T(H)=T,\; T_0(G)=T_0(H)=T_0,\; T_1(G)=T_1(H)=T_1;
\eqno(3.4)
$$
\noindent
in addition, $T_1\subseteq T_0$.

Let $G=\langle d_1,\text{Reg }G\rangle$ have the standard representation
$$
d_1=\sum_{\tau\in T_0}\dfrac{s_{\tau}}{m_{\tau}}e_0^{(\tau)},\eqno(3.5)
$$
where $\{e_0^{(\tau)}\,|\,\tau\in T_0\}$ is some $rc$-basis of the group $G$, $s_{\tau}$ is a $P_0(\tau)$-number, $\text{gcd}(s_{\tau},m_{\tau})=1$ ($\tau\in T_0$).

Let $H=\langle d_2,\text{Reg }H\rangle$ have the standard representation
$$
d_2=\sum_{\tau\in T_0}\dfrac{t_{\tau}}{m_{\tau}}f_0^{(\tau)},\eqno(3.6)
$$
where $\{f_0^{(\tau)}\,|\,\tau\in T_0\}$ is an $rc$-basis of the group $H$, $t_{\tau}$ is a $P_0(\tau)$-number, $\text{gcd}(t_{\tau},m_{\tau})=1$ ($\tau\in T_0$).

Let $s=(s_{\tau}\overline{1_{\tau}})_{\tau\in T_1}$, $t=(t_{\tau}\overline{1_{\tau}})_{\tau\in T_1}$.
It follow from \cite[Theorem 12.6,8]{Mad00} that $G\cong H$ if and only if $s=t\gamma v$ for some $\gamma\in\Gamma$ and $v\in V_{\infty}$. In other words, $G\cong H$ if and only if the elements $s$ and $t$ generate the same class $s\Gamma V_{\infty}=t\Gamma V_{\infty}$ in the quotient group $S^*/\Gamma V_{\infty}$.~$\square$

\textbf{Theorem 3.5.} Let $G,H\in\mathcal{A}_0$. The groups $G$ and $H$ are isomorphic if and only if the groups $\text{Mult}\,G$ and $\text{Mult}\,H$ are isomorphic.

\textbf{Proof.} If $G\cong H$, then it is clear that $\text{Mult}\,G\cong \text{Mult}\,H$.

Now let $\text{Mult}\,G\cong \text{Mult}\,H$. Then the groups $G$ and $H$ are near isomorphic by Theorem 3.3. According to Remark 3.2, we have $(3.4)$ for these groups.

Let the group $G=\langle d_1,\text{Reg}\,G\rangle$ have the standard representation $(3.5)$ and let the group $H=\langle d_2,\text{Reg}\,H\rangle$ have the standard representation $(3.6)$.

By Theorem 3.1 one of standard representations of the group $\text{Mult}\,G\cong \langle X,\text{Reg(Mult}\,G)\rangle$ is of the form
$$
X=\sum_{\tau\in T_0}\dfrac{s_{\tau}^{-1}}{m_{\tau}}F_0^{(\tau)},
$$
where $\{F_0^{(\tau)}\,|\,\tau\in T_0\}$ is an $rc$-basis of the group $\text{Mult}\,G$, $s_{\tau}^{-1}$ is a $P_0(\tau)$-number that is inverse to $s_{\tau}$ modulo $m_{\tau}$, $\tau\in T_0$.

Similarly, one of standard representations of the group 
$\text{Mult}\,H\cong \langle Y,\text{Reg(Mult}\,H)\rangle$
is of the form
$$
Y=\sum_{\tau\in T_0}\dfrac{t^{-1}}{m_{\tau}}J_0^{(\tau)},
$$
where $\{J_0^{(\tau)}\,|\,\tau\in T_0\}$ is an $rc$-basis of the group $\text{Mult}\,H$, $t_{\tau}^{-1}$ is a $P_0(\tau)$-number that is inverse to $t_{\tau}$ modulo $m_{\tau}$, $\tau\in T_0$.

Let $s=(s_{\tau}\overline{1_{\tau}})_{\tau\in T_1}$, $t=(t_{\tau}\overline{1_{\tau}})_{\tau\in T_1}$. In the group $S^*$, we have
$s^{-1}=(s_{\tau}^{-1}\overline{1_{\tau}})_{\tau\in T_1}$,
$t^{-1}=(t_{\tau}^{-1}\overline{1_{\tau}})_{\tau\in T_1}$.

Since $\text{Mult}\,G\cong \text{Mult}\,H$, it follows from Remark 3.4 that we have $s^{-1}=t^{-1}\gamma v$ in the group $S^*$ for some $\gamma\in\Gamma$ and $v\in V_{\infty}$. Consequently, $t=s\gamma v$.
According to Remark 3.4, the near isomorphic groups $G$ and $H$ are isomorphic.~$\square$

\section{Realization Theorem. $CRQ$-Groups Which Are Isomorphic to Their Multiplication Groups}\label{section4}

Our first aim is to describe groups that is implemented as the multiplication group of some group in $\mathcal{A}_0$.

\textbf{Remark 4.1.} Let $T$ be a set. We say that the system of positive integers $\{m_{\tau}\in \mathbb{N}\,|\,\tau\in T\}$ satisfies \textsf{condition $(m)$} if the following property is true: for any $p\in \mathbb{P}$, $k\in \mathbb{N}$ and $\tau\in T$, if $p^k$ divides $m_{\tau}$, then $p^k$ divides $m_{\sigma}$ for some $\sigma\in T\setminus\{\tau\}$.
We note that the system $\{m_{\tau}\,|\,\tau\in T\}$ satisfies condition $(m)$ if and only if the system $\{m_{\tau}\,|\,\tau\in T,\, m_{\tau}>1\}$ satisfies condition $(m)$.

Let $A=B\oplus C$, where $B$ and $C$ are two completely decomposable block-rigid groups of ring type and $B=\oplus_{\tau\in T(B)}R_{\tau}e_0^{(\tau)}$, where $R_{\tau}$ is a unitary subring of the field $\mathbb{Q}$, and the type of the additive group of $R_{\tau}$ is equal to $\tau$.

Let $\{m_{\tau}\,|\,\tau\in T(A)\}$ be a set of integers satisfying the following conditions:
\begin{itemize}
\item
$m_{\tau}$ is a $P_0$-number for any $\tau\in T(A)$;
\item
$m_{\tau}>1$ if and only if $\tau\in T(B)$.
\end{itemize}

In the group $\tilde{A}$, let we have the equality
$$
d=\sum_{\tau\in T(B)}\dfrac{s_{\tau}}{m_{\tau}}e_0^{(\tau)}\in \tilde{A}, \eqno (4.1)
$$
where $s_{\tau}$ is a $P _0(\tau)$-number that is co-prime to $m_{\tau}$ for every $\tau\in T(B)$. According to \cite[Remark 3.1, Remark 3.2]{KomT23}, the group $G=\langle d,A\rangle$ is a $CRQ$-group with regulator $A$, near isomorphism invariants $m_{\tau}$ $(\tau\in T(A))$ and the standard representation $(4.1)$ if and only if the set 
$\{m_{\tau}\,|\,\tau\in T(A)\}$ satisfies condition $(m)$.~$\square$

\textbf{Theorem 4.2.} For a group $M$, there exists a group $G\in \mathcal{A}_0$ such that $M\cong \text{Mult}\,G$ if and only if $M\in \mathcal{A}_0$ and for any $\tau\in T(M)$, we have 
$r(\text{Reg}_{\tau}M)=k_{\tau}^3$ for some $k_{\tau}\in \mathbb{N}$.

\textbf{Proof.} Let $M\cong \text{Mult}\,G$, where $G\in \mathcal{A}_0$. By Theorem 3.1, we have $M\in \mathcal{A}_0$, $T(M)=T(G)$, and $r(\text{Reg}_{\tau}M))=[r(\text{Reg}_{\tau}G)]^3$ for all $\tau\in T(M)$.

Let $M\in \mathcal{A}_0$, $T(M)=T$, $T_0(M)=T_0$, $m_{\tau}=m_{\tau}(M)$, $\tau\in T$. Let's assume that for every $\tau\in T$, we have $\text{Reg}_{\tau}M=\oplus_{\tau\in I_{\tau}}R_{\tau}e_i^{(\tau)}$, where
$$
I_{\tau}=\begin{cases}
\{0,\ldots,k_{\tau}^3-1\}, \text{ if } \tau\in T_0\\
\{1,\ldots,k_{\tau}^3\}, \text{ if } \tau\in T\setminus T_0.
\end{cases}, \quad k_{\tau}\in \mathbb{N},
$$ 
$R_{\tau}$ is a unitary subring of the field $\mathbb{Q}$, and the type of the additive group of $R_{\tau}$ is equal to $\tau$.

Let $E_0=\{e_0^{(\tau)}\,|\,\tau\in T_0\}$ be an $rc$-basis of the group $M$ and let the standard representation of the group $M$ be of the form 
$$
d=\sum_{\tau\in T_0}\dfrac{s_{\tau}}{m_{\tau}}e_0^{(\tau)},\eqno (4.2)
$$
where $s_{\tau}$ is a $P_0(\tau)$-number co-prime to $m_{\tau}$, $\tau\in T_0$. 
Then the system $\{m_{\tau}\,|\,\tau\in T_0\}$ satisfies condition $(m)$ by Remark 4.1.

We consider a completely decomposable group $A=\oplus_{\tau\in T}A_{\tau}$, where
$$
A_{\tau}=\oplus_{i\in J_{\tau}}R_{\tau}e_i^{(\tau)},\; J_{\tau}=\begin{cases}
\{0,\ldots,k_{\tau}-1\}, \text{ if } \tau\in T_0\\
\{1,\ldots,k_{\tau}\}, \text{ if } \tau\in T\setminus T_0.
\end{cases}
$$ 
For every $\tau\in T_0$, there exists a $P_0(\tau)$-number $s_{\tau}^{-1}$ that is inverse to $s_{\tau}$ modulo $m_{\tau}$. As such a number, one can take, for example, $s_{\tau}^{\varphi(m_{\tau})-1}$, where $\varphi(x)$ is the Euler function. We consider an element
$$
d_1=\sum_{\tau\in T_0}\dfrac{s_{\tau}^{-1}}{m_{\tau}}e_0^{(\tau)}\in \tilde{A}.\eqno (4.3)
$$
Since the system $\{m_{\tau}\,|\,\tau\in T_1\}$ satisfies condition $(m)$, it follows from Remark 4.1 that the group $G=\langle d_1,A\rangle$ is a group from $\mathcal{A}_0$ with regulator $A$, near isomorphism invariants $m_{\tau}$ ($\tau\in T(A)$) and the standard representation $(4.3)$.

It follows from Theorem 3.1 that for the group $\text{Mult}\,G$, we have 
$\text{Reg}(\text{Mult}\,G)\cong \text{Reg}\,M$ and 
$m_{\tau}(\text{Mult}\,G)=m_{\tau}(G)=m_{\tau}(M)$ for all $\tau\in T$. 
According to Remark 3.2, we obtain $M\cong_{nr}\text{Mult}\,G$.

By Theorem 3.1, the group $\text{Mult}\,G$ is isomorphic to the group $\langle X,M^{(2)}\rangle$ and one standard representations of the group $\langle X,M^{(2)}\rangle$ is of the form
$$
X=\sum_{\tau\in T_0}\dfrac{s_{\tau}}{m_{\tau}}F_0^{(\tau)}, \eqno (4.4) 
$$
where $\{F_0^{(\tau)}\,|\,\tau\in T_0\}$ is an $rc$-basis of the group $\text{Mult }G$. 

It follows from $(4.2)$ and $(4.4)$ that the near isomorphic groups $M$ and $\text{Mult}\,G$ are isomorphic by Remark 3.4.~$\square$

\textbf{Corollary 4.3.} For any rigid group $M\in \mathcal{A}_0$, there exists a group $G\in \mathcal{A}_0$ such that $M\cong \text{Mult}\,G$.~$\square$

We further will describe groups in the class $\mathcal{A}_0$ which are isomorphic to its multiplication group.

\textbf{Remark 4.4.} If the group $G\in\mathcal{A}_0$ is not rigid, then $r(G)\ne r(\text{Mult}\,G)$ by Theorem 3.1; therefore, $G$ cannot be not only isomorphic, but near isomorphic to the group $\text{Mult}\,G$. In this connection, we will further consider only rigid groups in the class $\mathcal{A}_0$. 

In the proofs we use the notation of Theorem 3.1 and Remark 3.4.

\textbf{Theorem 4.5.} If $G$ is a rigid $CRQ$-group of ring type, then the groups $G$ and $\text{Mult}\,G$ are near isomorphic.

\textbf{Proof.} Let $G$ be a rigid group in $\mathcal{A}_0$.
According to Theorem 3.1 we have $T(\text{Mult}\,G)=T(G)$, $m_{\tau}(\text{Mult}\,G)=m_{\tau}(G)$ for every $\tau\in T(G)$. Since $G$ is a rigid group, $r(\text{Reg}_{\tau}G)=1$ for every $\tau\in T(G)$. Consequently, $M_{\tau}^{(2)}=m_{\tau}^2A_{\tau}\cong A_{\tau}$ for all $\tau\in T(G)$. Therefore,
$$
\text{Reg}(\text{Mult}\,G)=M^{(2)}\cong A=\text{Reg}\,G.
$$
By Remark 3.2, the groups $G$ and $\text{Mult}\,G$ are near isomorphic.~$\square$

\textbf{Proposition 4.6.} Clipped direct summands in main decompositions of groups $G\in\mathcal{A}_0$ and $\text{Mult}\,G$ are near isomorphic.

\textbf{Proof.} Let $G\in\mathcal{A}_0$, $G_1$ be a clipped direct summand of a main decomposition of the group $G$ and let $\text{Reg}\,G_1=B=\oplus_{\tau\in T_0(G)}B_{\tau}$. Let $M_1$ be a clipped direct summand of a main decomposition of the group $\text{Mult}\,G$.

According to Theorem 3.1 and conditions $(2.1')$ and $(2.1'')$, we have
$$
T_0=T(G_1)=T_0(G)=T_0(\text{Mult}\,G)=T(M_1),
$$
$$
m_{\tau}=m_{\tau}(G_1)=m_{\tau}(G)=m_{\tau}(\text{Mult}\,G)=m_{\tau}(M_1)
$$
for all $\tau\in T_0$. In addition, it follows from \cite[Theorem 3.3]{KomT23} that
$$
\text{Reg}\,M_1\cong\oplus_{\tau\in T_0}
\begin{bmatrix}
m_{\tau}^2B_{\tau}&0&\ldots&0\\
0&0&\ldots&0\\
\ldots&\ldots&\ldots&\ldots\\
0&0&\ldots&0
\end{bmatrix}\cong \oplus_{\tau\in T_0}B_{\tau}=B=\text{Reg}\,G_1.
$$
Consequently, $G_1\cong_{nr} M_1$ by Remark 3.2.~$\square$

\textbf{Theorem 4.7.} Let $G$ be a group in $\mathcal{A}_0$ with standard representation
$$
d=\sum_{\tau\in T_0}\dfrac{s_{\tau}}{m_{\tau}}e_0^{(\tau)},
$$
where $T_0=T_0(G)$, $\{e_0^{(\tau)}\,|\,\tau\in T_0\}$ is an $rc$-basis of $G$, and let
$s=(s_{\tau}\overline{1_{\tau}})_{\tau\in T_0}\in S^*$. Then $G\cong \text{Mult}\,G$ if and only if $G$ is a rigid group and
$s^2\in \Gamma V_{\infty}$.

\textbf{Proof.} According to Remark 4.4, we consider a rigid group $G\in\mathcal{A}_0$. Then $T_1(G)=T_0(G)=T_0$. By Theorem 4.5, the groups $G$ and $\text{Mult}\,G$ are near isomorphic. By Theorem 3.1, we have $T_0(\text{Mult}\,G)=T_0(G)=T_0$ and there exists a standard representation of the group $\text{Mult}\,G$ which is of the form
$$
X=\sum_{\tau\in T_0}\dfrac{s_{\tau}^{-1}}{m_{\tau}}F_0^{(\tau)},
$$
where $\{F_0^{(\tau)}\,|\,\tau\in T_0\}$ is an $rc$-basis of the group $\text{Mult}\,G$, $s_{\tau}^{-1}$ is a $P_0(\tau)$-number that is inverse to $s_{\tau}$ modulo $m_{\tau}$, $\tau\in T_0$. 
Then we have
$$
s^{-1}=(s_{\tau}^{-1}\overline{1_{\tau}})_{\tau\in T_0}
$$
in the group $S^*$.

According to Remark 3.4, for near isomorphic groups $G$ and $\text{Mult}\,G$, the isomorphism $G\cong\text{Mult}\,G$ holds if and only if the equality $s^{-1}\Gamma V_{\infty}=s\Gamma V_{\infty}$ holds in the quotient group $S^* /\Gamma V_{\infty}$. This condition is equivalent to the inclusion $s^2\in \Gamma V_{\infty}$.~$\square$

So-called \textsf{proper} groups provide an example of groups in the class $\mathcal{A}_0$ that are isomorphic to their multiplication groups. According to \cite{BlaIS01}, a group $G\in \mathcal{A}_0$ is said to be \textsf{proper} if it admits the standard representation $(2.3)$ such that $s_{\tau}=1$ for all $\tau\in T_0(G)$. In the class $\mathcal{A}_0$, proper groups play a special role. For example, for any group $G\in \mathcal{A}_0$, the endomorphism group $\text{End}\,G$ is a proper group from $\mathcal{A}_0$.

\textbf{Corollary 4.8.} If $G$ is a rigid proper group in the class $\mathcal{A}_0$, then the groups $G$ and $\text{Mult}\,G$ are isomorphic.~$\square$

If $G\in \mathcal{A}_0$ and $r(G)=1$, then $G\cong \text{Mult}\,G$, e.g., see \cite{KomN23}.
We will show that for any $k\in \mathbb{N}$, $k\ge 2$, there exists a rigid group $G$ of rank $k$ in the class $\mathcal{A}_0$ that is not isomorphic to the group $\text{Mult}\,G$.

\textbf{Example 4.9.} Let $k\in \mathbb{N}$, $k\ge 2$, and let $p$ be a prime number, $p>3$. By Dirichlet's theorem \cite[Theorem 3.3.1]{FinR16} the arithmetic progression $\{1+pt\,|\,t\in \mathbb{N}\}$ contains infinitely many prime numbers, so we choose distinct numbers $q_1,q_2,\ldots,q_k,s_1$ in this progression.
Let $s_2$ be an integer such that $1<s_2<p-1$. Since $s_2\not\equiv 1\,(\text{mod }p)$ and $s_2\not\equiv -1\,(\text{mod }p)$, we have
$$
s_2^2\not\equiv 1\, (\text{mod }p)\eqno (4.5)
$$
by \cite[Lemma 2.5.2.1]{FinR16}. In addition, $s_2<p-1<q_2$; therefore $\text{gcd}(s_2,q_2)=1$. We choose positive integers $s_3,\ldots,s_k$ such that $\text{gcd}(s_i,q_i)=\text{gcd}(s_i,p)=1$, $i=3,\ldots,k$.

We consider idempotent types $\tau_1,\ldots, \tau_k$ such that $P_{\infty}(\tau_i)=\{q_i\}$, ($i=1,\ldots,k$). Let $A=\oplus_{i=1}^kR_ie_i$, where $R_i$ is a unitary subring of the field $\mathbb{Q}$ with additive group of type $\tau_i$. Then $A$ is a rigid group of ring type.

Let $m_1=\ldots =m_k=p$. Then $s_i$ and $m_i$ are co-prime $P_0(\tau_i)$-numbers for every $i\in \{1,\ldots,k\}$. In the group $\tilde{A}$, we consider an element
$$
d=\sum_{i=1}^k\dfrac{s_i}{m_i}e_i.\eqno (4.6)
$$
Since the system $\{m_1,\ldots,m_k\}$ satisfies condition $(m)$, it follows from Remark 4.1 that $G=\langle d,A\rangle$ is a group in $\mathcal{A}_0$ with regulator $A$, near isomorphism invariants $m_1,\ldots,m_k$, and standard representation $(4.6)$.

Let's assume that $G\cong\text{Mult}\,G$. By Theorem 4.7, we have $s^2=\gamma v$, where $s=(s_i\overline{1_i})_{i=1,\ldots,k}\in S^*$, $\gamma\in\Gamma$, $v\in V_{\infty}$.

Since $m_1=\ldots =m_k=p$, the element $\gamma$ can be represented in the form 
$\gamma=(\alpha\overline{1_1},\ldots,\alpha\overline{1_k})$ for some $\alpha\in\mathbb{Z}$. In addition, since $P_{\infty}(\tau_i)=\{q_i\}$, we have $v=(q_1^{t_1}\overline{1_1},\ldots,q_k^{t_k}\overline{1_k})$, where 
$t_i\in \mathbb{Z}$ and $q_i^{-1}$ is an integer that is inverse to $q_i$ modulo $p$. Consequently, 
$$
s_1^2\equiv \alpha q_1^{t_1} (\text{mod }p),\;
s_2^2\equiv \alpha q_2^{t_2} (\text{mod }p).
$$
Since $q_1\equiv 1 (\text{ mod }p)$ and $q_2\equiv 1 (\text{ mod }p)$, we have $s_1^2\equiv s_2^2 (\text{mod }p)$. This contradicts to $(4.5)$, since $s_1\equiv 1\, (\text{mod }p)$. Consequently, the groups 
$G$ and $\text{Mult}\,G$ are not isomorphic.~$\square$

\end{document}